\documentclass[12pt,draftcls,onecolumn,web]{ieeecolor}
\usepackage{generic}
\usepackage{cite}
\usepackage{amsmath,amssymb,amsfonts}
\usepackage{algorithmic}
\usepackage{graphicx}
\usepackage{textcomp}
\usepackage{array}
\usepackage{caption}
\usepackage{color}
\allowdisplaybreaks[4]

\def\BibTeX{{\rm B\kern-.05em{\sc i\kern-.025em b}\kern-.08em
    T\kern-.1667em\lower.7ex\hbox{E}\kern-.125emX}}
\begin{document}
\title{Distributed Control of Linear Quadratic Mean Field Social Systems with Heterogeneous Agents}
\author{Yong Liang, Bing-Chang Wang, \IEEEmembership{Senior Member, IEEE}, and Huanshui Zhang, \IEEEmembership{Senior Member, IEEE}
\thanks{This work was supported by the National Natural Science Foundation of China under Grants 61922051, 62122043, 62192753, 
	Major Basic Research of Natural Science Foundation of Shandong Province (ZR2021ZD14, ZR2020ZD24), 
	Science and Technology Project of Qingdao West Coast New Area (2019-32, 2020-20, 2020-1-4), High-level Talent Team Project of Qingdao West Coast New Area (RCTD-JC-2019-05), 
	and Key Research and Development Program of Shandong Province (2020CXGC01208).}
\thanks{Y. Liang is with the School of Control Science and Engineering,
	Shandong University, 	Jinan, China (e-mail: yongliang@mail.sdu.edu.cn). }
\thanks{B.-C. Wang is with the School of Control Science and Engineering,
	Shandong University, 	Jinan, China (e-mail: bcwang@sdu.edu.cn).}
\thanks{H. Zhang is with 
the College of Electrical Engineering and Automation, Shandong University of Science and Technology, Qingdao, China (e-mail: hszhang@sdu.edu.cn).}}

\maketitle

\newtheorem{problem}{Problem}
\newtheorem{lemma}{Lemma}
\newtheorem{proposition}{Proposition}
\newtheorem{theorem}{Theorem}

\begin{abstract}
In this paper, we study the social optimality for mean field linear-quadratic control systems following the direct approach, where subsystems are coupled via individual dynamics and costs according to a network topology. A graph is introduced to represent the network topology of the large-population system, where nodes represent subpopulations called clusters and edges represent communication relationship. By the direct approach, we first seek the optimal controller under centralized information structure, which characterized by a set of forward-backward stochastic differential equations. Then the feedback controller is obtained with the help of Riccat equations. Finally, we design the distributed controller with mean field approximations, which has the property of asymptotically social optimality.

\end{abstract}

\section{Introduction}
In recent years, mean field games and control have became a hot topic, and it has wide applications in system control, applied mathematics, and economics \cite{b1,b2,b3}, such as dynamic production adjustment \cite{b4}, vaccination games \cite{b5}, resource allocation in internet of things \cite{b6}, etc. Mean field games theory originated from the parallel works of Lasry and Lions \cite{b7} and of Huang et al. \cite{b8}. A key feature of the mean field system is its mean field coupling terms, which often appears in multi-agent systems, such as the team decision problem with partially exchangeable agents \cite{b9}. 

Under the centralized information structure, the computational complexity increases exponentially with the increase of the number of agents. Due to communication constraints, it is also unrealistic for each agent to obtain global information. Therefore, it is meaningful to design decentralized controllers only using local information. The difficulty is how to use local information to approximate the mean field coupling term which is related to global information. At present, there are mainly two approaches: fixed point and direct approaches. The first approach uses an infinite to finite population approach. First the fixed point equation is obtained in the infinite population system and then substitute the solution to the finite population system to analyze the property of $\varepsilon$-Nash equilibrium. See \cite{b7,b10,b11,b12,b16} for details. The second approach uses an finite to infinite population approach. First, one directly solve the Nash equilibrium under the finite population and then obtain the limit of the  Nash equilibrium with respect to the number of populations tending to infinity. Finally, it is proved that the limit is a $\varepsilon$-Nash equilibrium in the original system. See \cite{b8,b13,b15,b16,b17} for details. As the complexity of the system increases, the fixed point equation will become more complicated and difficult to solve, such as mean field games with major-minor agents \cite{b16}. Even the fixed point approach will fail when the system contains common noise \cite{b17}.

In this paper, we study the mean field social control problem by the direct approach. The agents in the large-population system are heterogeneous, which makes it difficult to solve the centralized controller. The case we consider in this paper is that the large-population system contains $K$-type agents. Agents of the same type are divided into a cluster and represented as a node of the graph $\mathcal{G}_K$. The agents in the same cluster generate a local mean field term. Therefore there are total of $K$ local mean field terms in this large-population system. Agents are coupled with each other in dynamics and costs through the mean field terms according to the adjacency matrix. Compared with the mean field social control involving homogeneous agents, the heterogeneous case is difficult to analyze and solve the centralized controller only using traditional methods. The key to solve this problem is to construct the local mean field terms. The distributed controller is designed by the optimal estimation of $K$ local mean field terms. The main contributions of the paper are listed as follows.
\begin{itemize}
	\item By variational analysis, the necessary and sufficient condition for the exist for the existence of centralized controllers is obtained, which is characterized by the adapted solution of a set of forward-backward stochastic differential equations.
	\item The centralized optimal feedback controller is obtained by introducing Riccati equations.  
	\item A distributed controller is proposed by mean field approximations, which has the property of asymptotically social optimality.   
	
\end{itemize}

The rest of the paper is organized as follows. In section \ref{sec1} we formulate the mean field social control problem. Section \ref{centralized} gives the centralized results of this paper. Section \ref{distributed} design distributed controllers based on centralized results and show the asymptotically social optimality. Section \ref{sec3} concludes this paper.  

Notations: For a matrix $M$, $M^T$ means the transpose of $M$. For a set of vectors $\{x_1, x_2, . . .,x_N\}$, $\mathrm{vec}(x_1, x_2, . . .,x_N)$ denotes the vector $(x_1^T,x_2^T, . . .,x_N^T)^T$. For a set of matrices $\{A_1,A_2, . . .,A_N\}$, $\mathrm{rows}(A_1,$ $A_2, . . .,A_N)$, $\mathrm{cols}(A_1,A_2, . . .,A_N)$, and $\mathrm{diag}(A_1,A_2, . . .,A_N)$ denote the matrices $(A_1^T,A_2^T, . . .,A_N^T)^T$, $(A_1,$ $A_2, . . .,A_N)$ and the diagonal matrix with the elements of $(A_1,A_2, . . .,A_N)$ on the main diagonal.  Let $I_n$ be the $n$-dimensional identity matrix. Let $1_{n\times m}$ and $0_{n\times m}$ denote the $n\times m$ matrix with all the elements $1$ and $0$, respectively. We use $C_0,C_1,C_2$ etc. to denote generic constants, which may vary from place to place. For two matrices $A=[a_{ij}]$ and $B=[b_{ij}]$, let $\otimes$ and $\odot$ denote the Kronecker product and the dot product, respectively, i.e.,
\begin{equation*}
\begin{split}
		A\otimes B=& \left(\begin{array}{cccc}
		a_{11}B&. . .&a_{1n}B\\. . .&. . .&. . .\\a_{m1}B&. . .&a_{mn}B
	\end{array} \right),
\quad
A\odot B={\small \left(\begin{array}{cccc}
	a_{11}b_{11}&. . .&a_{1n}b_{1n}\\. . .&. . .&. . .\\a_{m1}b_{m1}&. . .&a_{mn}b_{mn}
\end{array} \right)}.
\end{split}
\end{equation*}

\section{Problem Formulation}\label{sec1}

\subsection{Network topology}


Consider a large-population system with $K\in\mathbb{R}$ clusters $\mathcal{C}_q$, $q\in\mathcal{V}_K\triangleq \{1,2, . . .,K\}$, where the cluster $C_q$ contains $N_q$ homogeneous agents with the same dynamics, cost functions and communication capabilities. In this scenario, the large-population system contains a total of $N\triangleq \sum_{q\in\mathcal{V}_K}N_q$ heterogeneous agents $\mathcal{A}_i$, $i\in \mathcal{I}\triangleq \{1,2, . . .,N\}$, and the agents are partitioned into the $K$ disjoint clusters $\mathcal{C}_q$, $q\in\mathcal{V}_K$. The set of agents belonging to the cluster $\mathcal{C}_q$ is denoted by $\mathcal{I}_q\subseteq \mathcal{I}$ with $|\mathcal{I}_q|=N_q$, $\cup_{q\in\mathcal{V}_K}\mathcal{I}_q=\mathcal{I}$, and $\cap_{q\in\mathcal{V}_K}\mathcal{I}_q=\emptyset$. In this paper, we also list the agents as $\mathcal{A}_i$, $i\in\mathcal{I}_q$, $q\in\mathcal{V}_K$.

The agents in different clusters communicate with each other through a network whose topology is modeled as a directed graph $\mathcal{G}_K=\{ \mathcal{V}_K,\mathcal{E}_K,\mathcal{M}_K\}$, where the clusters and the communication channels between clusters are represented by the node set $\mathcal{V}_K$ and the edge set $\mathcal{E}_K\subseteq\mathcal{V}_K\times\mathcal{V}_K$, respectively.  An edge denoted by the pair $(j,i)$ represents a  communication channel from cluster $\mathcal{C}_j$ to cluster $\mathcal{C}_i$. The neighbor set of cluster $\mathcal{C}_q$ is denoted by $\mathcal{N}_q\triangleq \{\ p\  |\ p\in \mathcal{V}_K,\ (p,q)\in \mathcal{E}_K\}$.  Let $E_K=[e_{ij}]\in\mathbb{R}^{K\times K}$ be the communication matrix, where $e_{ij}=1$ if $(j,i)\in \mathcal{E}_K$, $e_{ij}=0$ if $(j,i)\notin \mathcal{E}_K$. The weighted adjacency matrix of $\mathcal{G}_K$ is denoted by $\mathcal{M}_K=[m_{ij}]\in\mathbb{R}^{K\times K}$.  The agents are coupled with each other according to the adjacency matrix. 

\subsection{Dynamics and costs}
Let $x_i\in\mathbb{R}^n$ and $u_i\in\mathbb{R}^m$ be the state and control of agent  $\mathcal{A}_i,i\in\mathcal{I}_q, q\in\mathcal{V}_K  $. The local cluster mean field term of cluster $\mathcal{C}_q$ is defined as the  empirical mean of the states of all agents in the cluster, i.e.,
\begin{align}
	x^K_q(t)\triangleq \frac{1}{N_q}\sum_{i\in\mathcal{I}_q}x_i(t),\quad q\in\mathcal{V}_K.
\end{align}
Therefore the global cluster mean field term of the large-population system is stacked as $x^K\triangleq \mathrm{vec}(x^K_1,x^K_2,$ $ . . .,x^K_K)$. By the adjacency matrix $\mathcal{M}_{K}$, the influence of the global cluster mean field term $x^K$ on the agents in cluster $C_q$ is denoted as $z^K_q$, where 
\begin{align}
	z^K_q(t)\triangleq \frac{1}{K}\sum_{p\in\mathcal{V}_{K}}m_{qp}x^K_p(t),\quad q\in\mathcal{V}_K.
\end{align}
Coupled by the term $z^K_q$, the dynamics of agent $\mathcal{A}_i,i\in\mathcal{I}_q,q\in\mathcal{V}_q$, is given by the following stochastic differential equation
\begin{equation}\label{e122701}
		\begin{aligned}dx_i(t)=&\big[A_qx_i(t)+B_qu_i(t)+G_qz_q^K(t) \big]dt+\Sigma_q dw_i(t),\quad x_i(0)=x_{i0},
	\end{aligned}
\end{equation}
where $A_q, B_q, G_q, \Sigma_q$ are constant matrices with compatible dimensions and $\{w_i,i\in\mathcal{I}_q\}$ are $1$-dimensional independent standard Brownian motions defined on a complete filtered probability space $(\Omega,\mathcal{F},\mathbb{P})$. 

Let $\mathbf{x}\triangleq \mathrm{vec}(x_1,x_2, . . .,x_N)$ and $\mathbf{u}\triangleq \mathrm{vec}(u_1,u_2, . . .,u_N)$ be the state and control of the large-population system, respectively. Then the cost function of agent $\mathcal{A}_i, i\in\mathcal{I}_q, q\in\mathcal{V}_K$, is given by
\begin{equation}\label{e122802}
	\begin{split}
			J_i(\mathbf{u})=&\mathbb{E}\int_0^T \big[|x_i(t)-\Gamma_qz_q^K(t)|^2_{Q_{q}}+|u_i(t)|^2_{R_{q}}\big]dt+\mathbb{E}|x_i(T)-\Gamma_qz_q^K(T)|^2_{H_q}, \nonumber 
	\end{split}
\end{equation}
where $\Gamma_q,Q_q\geq0,R_q>0,H_q\geq 0$ are constant matrices with compatible dimensions. The social cost function is defined as 
\begin{align}
	&J_{\rm soc}^{\rm (N)}(\mathbf{u})=\sum_{q\in\mathcal{V}_K}\sum_{i\in\mathcal{I}_q}J_i(\mathbf{u}).\label{e122702} 
\end{align}
Let $M_q$, $q\in\mathcal{V}_K$, be the $q$-th row of $\mathcal{M}_K$ and 
\begin{equation*}
	\bar{G}_q\triangleq \frac{M_q}{K}\otimes G_q.
\end{equation*}
Then the individual dynamics \eqref{e122701} can be rewritten as 
\begin{equation}\label{e122801}
	\begin{split}
			dx_i(t)=&\big[A_qx_i(t)+B_qu_i(t)+\bar{G}_qx^K(t) \big]dt+\Sigma_q dw_i(t),\quad x_i(0)=x_{i0}.
	\end{split}
\end{equation}
Similarly, let 
\begin{equation*}\label{e072002}
\bar{\Gamma}_q\triangleq \frac{M_q}{K}\otimes \Gamma_q.
\end{equation*}
Thus the social cost \eqref{e122702} can be rewritten as 
\begin{equation}\label{e122803}
	\begin{split}
			J_{\rm soc}^{\rm (N)}(\mathbf{u})=&\sum_{q\in\mathcal{V}_K}\sum_{i\in\mathcal{I}_q}\Big\{\mathbb{E}\int_0^T \big[|x_i(t)-\bar{\Gamma}_qx^K(t)|^2_{Q_{q}}+|u_i(t)|^2_{R_{q}}\big]dt+\mathbb{E}|x_i(T)-\bar{\Gamma}_qx^K(T)|^2_{H_q}\Big\}. 
	\end{split}
\end{equation}

\subsection{Main problems}
In this paper, we apply the direct approach to sequentially study the mean field social control problem under centralized and distributed information patterns. So we first give the definitions of filtration used in this paper. Denote 
\begin{align*}
	\mathcal{F}_t^i\triangleq& \sigma( x_i(s),w_i(s),s\leq t),\\
	\mathcal{F}_t\triangleq& \sigma( \cup_{i\in\mathcal{I}_q,q\in\mathcal{V}_K} \mathcal{F}_t^i),\\
	\mathbb{F}\triangleq &\{\mathcal{F}_t\}_{0\leq t\leq T},\\
	\mathcal{H}_t^{q,K}\triangleq& \sigma(\cup_{j\in\mathcal{I}_p,p\in\mathcal{N}_q} \mathcal{F}^j_t),\\
	\mathcal{H}_t^i\triangleq &\sigma(\mathcal{H}_t^{q,K},\mathcal{F}_t^i),\\
	\mathbb{H}_i\triangleq& \{ \mathcal{H}^i_t\}_{0\leq t\leq T}.
\end{align*}
Then the centralized and distributed admissible control sets of agent $\mathcal{A}_i,i\in\mathcal{I}_q,q\in\mathcal{V}_K$, are given by $\mathcal{U}_c[0,T]$ and $\mathcal{U}_{d}^i[0,T]$, respectively, where
\begin{align*}
	\mathcal{U}_{c}[0,T]=&\big\{u(\cdot)\in L^2_{\mathbb{F}}(0,T;\mathbb{R}^{n})\big\},\\	\mathcal{U}_{d}^i[0,T]=&\big\{  u(\cdot)\in L^2_{\mathbb{H}_i}(0,T;\mathbb{R}^{n})   \big\}.
\end{align*}
As stated in the Introduction, to our knowledge, there is currently no general method for designing optimal distributed controllers. Thus we consider asymptotically optimal distributed  controllers based on the mean field approximation methodology. To characterize the asymptotic optimality of distributed controllers, the following definition is introduced: 

\noindent\textbf{Definition 1:} For the dynamics \eqref{e122701} with the social cost function \eqref{e122702}, a set of control $\check{u}_i\in\mathcal{U}_d^i[0,T], i\in\mathcal{I}_q,q\in\mathcal{V}_K$, is called an asymptotically optimal distributed controller with respect to the number of agents in the large-population system if 
\begin{align}
	\Big|\frac{1}{N}J_{\rm soc}^{\rm (N)}(\check{\mathbf{u}})-\frac{1}{N}\inf_{u_i\in\mathcal{U}_c[0,T]}J_{\rm soc}^{\rm (N)}(\mathbf{u})\Big|=o(1).
\end{align}
We propose the following two problems:  
\begin{problem}[P1]\label{p1}
	For each agent $\mathcal{A}_i$, $i\in \mathcal{I}_q, q\in\mathcal{V}_K$, 
	  find a $\mathcal{F}_t$-adapted  optimal centralized controller $\hat{u}_i\in \mathcal{U}_{c}[0,T]$ to minimize the social cost function \eqref{e122702}.
\end{problem}
\begin{problem}[P2]\label{p2}
	For each agent $\mathcal{A}_i$, $i\in \mathcal{I}_q, q\in\mathcal{V}_K$, 
	  find a $\mathcal{H}_t^i$-adapted asymptotically optimal  distributed controller $\check{u}_i\in \mathcal{U}_{d}^i[0,T]$ to minimize the social cost function \eqref{e122702}.
\end{problem}

Therefore, in this paper we first solve Problem (P1) to obtain  optimal centralized controllers, and then further study Problem (P2) to design distributed asymptotically optimal controllers based on mean field approximations. The following assumption on the distribution of initial states is imposed.

(A1) The initial states $x_{i0}$, $i\in\mathcal{I}_q$, $q\in\mathcal{V}_K$, are mutually independent with $\mathbb{E}[x_i(0)]=\bar{m}_q$, $i\in\mathcal{I}_q$, $q\in\mathcal{V}_K$, and there exists a finite constant $C_0$ independent of $N$ such that $\max_{1\leq i\leq N}\mathbb{E}|x_i(0)|^2\leq C_0$. Let $\bar{m}^K\triangleq\mathrm{diag}(\bar{m}_1, \bar{m}_2, . . .,\bar{m}_K)$.



\section{Optimal Centralized  Controllers}\label{centralized}

In this section we first solve (P1) by variational analysis to obtain the open-loop centralized optimal controller which is characterized by a set of FBSDEs, then obtain its feedback representation by virtue of two Riccati equations. 
\subsection{Open-loop controllers and MF-FBSDEs}
We make some denotations for convenience of the derivation of the main result. Denote $\pi_q\triangleq N_q/N$, $q\in\mathcal{V}_K$. Then $\pi=(\pi_1,\pi_2, . . .,\pi_K)$ is a probability vector which gives the empirical distribution of the $K$-type agents. Denote
\begin{align*}
	\Pi^K\triangleq&\ \mathrm{diag}(I_n\otimes\pi_1,I_n\otimes\pi_2, . . .,I_n\otimes\pi_K),\\
	N^K\triangleq&\ \mathrm{diag}(I_{n}\otimes N_1,I_{n}\otimes N_2, . . .,I_{n}\otimes N_K),\\
	G^K\triangleq &\ \mathrm{rows}(\bar{G}_1,\bar{G}_2, . . .,\bar{G}_K),\\ Q^K\triangleq& \mathrm{diag}(Q_1,Q_2, . . .,Q_K),\\
	\Gamma^K\triangleq&\ \mathrm{rows}(\bar{\Gamma}_1,\bar{\Gamma}_2, . . .,\bar{\Gamma}_K),\\ H^K\triangleq& \mathrm{diag}(H_1,H_2,. . .,H_K),\\
	D^K\triangleq&\ \mathrm{cols}(D_1,D_2, . . .,D_K)\\=&\ N^KG^K(N^K)^{-1},\\
\bar{Q}^K\triangleq &\ \mathrm{rows}(\bar{Q}_1,\bar{Q}_2, . . .,\bar{Q}_K)\\=&\   Q^K\Gamma^K+(N^K)^{-1}(\Gamma^K)^TQ^KN^K -(N^K)^{-1}(\Gamma^K)^TN^KQ^K\Gamma^K,\\
\bar{H}^K\triangleq &\ \mathrm{rows}(\bar{H}_1,\bar{H}_2, . . .,\bar{H}_K)\\=&\   H^K\Gamma^K+(N^K)^{-1}(\Gamma^K)^TH^KN^K -(N^K)^{-1}(\Gamma^K)^TN^KH^K\Gamma^K.
  \end{align*}
Note that we have $N^K=N \Pi^K$. By variational analysis, the necessary and sufficient condition is obtained for the solvability of Problem (P1) as follows.
\begin{theorem}
	(P1) is solvable if and only if the following MF-FBSDEs 
	\begin{align}
		\begin{split}
			\left\{ 
			\begin{aligned}
				&d\hat{x}_i(t)=\big[A_q\hat{x}_i(t)+B_q\hat{u}_i(t)+\bar{G}_q\hat{x}^K(t) \big]dt+\Sigma_q dw_i(t),\\
				&d\lambda_{i}(t)= -[A_q^T\lambda_i(t)+D_q^T\lambda^K(t)+ Q_q\hat{x}_i(t)-\bar{Q}_q\hat{x}^{K}(t)] dt+\sum_{q\in\mathcal{V}_K}\sum_{j\in\mathcal{I}_q}\beta_{i}^j(t)dw_j(t),\\& \hat{x}_i(0)=x_{i0}, \ \lambda_i(T)=H_q\hat{x}_i(T)-\bar{H}_q\hat{x}^{K}(T), i\in\mathcal{I}_q,q\in\mathcal{V}_K\label{e122804}
			\end{aligned}
			\right.
		\end{split}
	\end{align}
	admit a $\mathcal{F}_t$-adapted solution $(\hat{x}_i(t),\lambda_i(t),i\in\mathcal{I}_q,q\in\mathcal{V}_K)$ on $t\in[0,T]$, where $\lambda^K\triangleq \mathrm{vec}(\lambda^K_1,\lambda^K_2, . . .,\lambda^K_K)$ with $\lambda^K_q\triangleq (1/N_q)\sum_{i\in\mathcal{I}_q}\lambda_i$. The optimal centralized open-loop controller for agent $\mathcal{A}_i,i\in\mathcal{I}_q,q\in\mathcal{V}_K$, is given by 
	\begin{align}
		\hat{u}_i(t)=-R_q^{-1}B_q^T\lambda_i(t).\label{e122901} 
	\end{align}

\end{theorem} 

\subsection{Feedback representation and Riccati equations}
Based on the open-loop centralized optimal controller \eqref{e122901}, we next obtain its feedback representation by introducing two Riccati equations. Moreover, we also make the following denotations for convenience of discussions. 
\begin{align*}
	A^K&\triangleq \mathrm{diag}(A_1,A_2, . . .,A_K),\\
	B^K&\triangleq \mathrm{diag}(B_1,B_2,. . .,B_K),\\
	R^K&\triangleq\mathrm{diag}(R_1,R_2, . . .,R_K),\\
	\Sigma^K&\triangleq \mathrm{diag}(\Sigma_1,\Sigma_2, . . .,\Sigma_K).
\end{align*}
We propose the following Riccati equations 
\begin{subequations}
	\begin{align}
	\begin{split}\label{e122902}
			0=&\dot{P}^K(t)+P^K(t)A^K+(A^K)^TP^K(t)+Q^K-P^K(t)B^K(R^K)^{-1}(B^K)^TP^K(t),
	\end{split}\\
		\begin{split}\label{e031601}
				0=&\dot{K}^K(t)+[A^K+N^KG^K(N^K)^{-1}]^TK^K(t)+K^K(t)[A^K+G^K]+P^K(t)G^K\cr&+(N^K)^{-1}(G^K)^TN^KP^K(t)-\bar{Q}^K-K^K(t)B^K(R^K)^{-1}(B^K)^T[P^K(t)+K^K(t)]\cr&-P^K(t)B^K(R^K)^{-1}(B^K)^TK^K(t),
		\end{split}\\
			 P^K&(T)=H^K,\quad K^K(T)=-  \bar{H}^K.\nonumber
	\end{align}	
\end{subequations}
Let $P_q$ be the $q$-th diagonal element of $P^K$ and $\bar{K}_q$ be the $q$-th row of $K^K$, $q\in\mathcal{V}_K$. Then the Riccati equations  \eqref{e122902}-\eqref{e031601} can be rewritten as 
	\begin{subequations}
	\begin{align}
	\begin{split}\label{e122904}
			0=&\dot{P}_q(t)+P_q(t)A_q+A_q^TP_q(t)+Q_q-P_q(t)B_qR_q^{-1}B_q^TP_q(t),
	\end{split}\\
	\begin{split}\label{e122905}
			0=&\dot{\bar{K}}_q(t)+A_q^T\bar{K}_q(t)+D_q^TK^K(t)+\bar{K}_q(t)(A^K+G^K)+P_q(t)\bar{G}_q+D_q^TP^K(t)\cr&-\bar{K}_q(t)B^K(R^K)^{-1}(B^K)^T[P^K(t)+K^K(t)]-\bar{Q}_q-P_q(t)B_qR_q^{-1}B_q^T\bar{K}_q(t),
	\end{split}\\
	 P_q&(T)=H_q,\quad \bar{K}_q(T)=-\bar{H}_q,\quad q\in\mathcal{V}_K.	\nonumber
	\end{align}
\end{subequations}
Therefore the result of feedback representation is given as follows.
\begin{theorem}
	(P1) admits optimal centralized feedback controllers if and only if the Riccati equations \eqref{e122902}-\eqref{e031601} admit solutions $P^K(t),K^K(t)$ on $t\in[0,T]$. In this case, the feedback controller of agent $\mathcal{A}_i$, $i\in\mathcal{I}_q$, $q\in\mathcal{V}_K$, is given by 
	\begin{align}
		\hat{u}_i(t)=&-R_q^{-1}B_q^TP_q(t)\hat{x}_i(t)-R_q^{-1}B_q^T\bar{K}_q(t)\hat{x}^{K}(t),\label{e122909}
	\end{align}
	where $P_q$ and $\bar{K}_q$ are given by \eqref{e122904} and \eqref{e122905}. The value function admits the following form:
	\begin{equation}\label{e032305}
		{\rm V}(\mathbf{x}_0)= \sum_{q\in\mathcal{V}_K}\sum_{i\in\mathcal{I}_q}\mathbb{E}\langle P_q(0)x_{i0},x_{i0}\rangle + \mathbb{E}\langle N^KK^K(0)x^K_0,x^K_0\rangle.
	\end{equation}
	
\end{theorem}

We introduce the following notations to simplify the expression of closed-loop systems.
\begin{align}
\tilde{A}_q\triangleq&\ A_q-B_qR_q^{-1}B_q^TP_q,\cr
\tilde{G}_q\triangleq&\ \bar{G}_q-B_qR_q^{-1}B_q^T\bar{K}_q,\cr
\tilde{A}^K\triangleq&\ \mathrm{diag}(\tilde{A}_1,\tilde{A}_2,. . . ,\tilde{A}_K),\cr 
\tilde{G}^K\triangleq&\ \mathrm{rows}(\tilde{G}_1,\tilde{G}_2, . . .,\tilde{G}_K).\nonumber
\end{align}
Under the feedback controller \eqref{e122909}, the closed-loop systems of $\hat{x}_i,\hat{x}^K_q$ and $\hat{x}^K$, $i\in\mathcal{I}_q,q\in\mathcal{V}_K$ are given by
\begin{subequations}\label{e010301}
	\begin{align}	
		\begin{split}\label{e010302}
			d\hat{x}_i(t)=&\ 
			 [\tilde{A}_q(t)\hat{x}_i(t)+\tilde{G}_q(t)\hat{x}^{K}(t)]dt+\Sigma_q dw_i(t),\quad \hat{x}_i(0)=x_{i0},
		\end{split}\\
		\begin{split}\label{e010303}
			d\hat{x}_q^K(t)=&\ [\tilde{A}_q(t)\hat{x}_q^K(t)+\tilde{G}_q(t)\hat{x}^{K}(t) ]dt+\Sigma_q dw_q^K(t),\quad \hat{x}_q^K(0)=x^K_{q0},
		\end{split}\\
		\begin{split}\label{e010304}
			d\hat{x}^K(t)=&\ [\tilde{A}^K(t)+\tilde{G}^K(t)]\hat{x}^{K}(t)dt+\Sigma^K dw^K(t),\quad \hat{x}^K(0)=x^K_0. 
		\end{split} 
	\end{align}
\end{subequations}

\section{Asymptotically Optimal Distributed Controllers}\label{distributed}
In this section we will design the asymptotically optimal distributed controller based on the optimal centralized feedback controller \eqref{e122909} according to the network topology $\mathcal{G}_K$. We first design cluster mean field estimators to estimate the global cluster mean filed term under the distributed information pattern. Then distributed controllers are proposed by using the cluster mean field estimators. Finally, we prove the asymptotic optimality of the distributed controllers. We use $\check{x}_i$ to denote the state of agent $\mathcal{A}_i$, $i\in\mathcal{I}_q$, $q\in\mathcal{V}_K$, under the distributed controller $\check{u}_i$. Let  $\check{x}^K_q$ and $\check{x}^K$ be the corresponding local and global cluster mean field terms, respectively.
\subsection{Cluster mean field estimators}
Note that the optimal centralized feedback controller \eqref{e122909} contains the individual state and the global cluster mean field term. Moreover, the ability of each agent to acquire information is different according to the communication matrix $E_K$. Therefore, the main idea of the design of distributed controllers is that each cluster makes a local estimation of the global cluster mean field term according to the network topology. 

Therefore we design cluster mean filed estimators for agents in each cluster to estimate the global cluster mean filed term. In the following context we denote the local estimation of agent $\mathcal{A}_i$ in cluster $\mathcal{C}_q,q\in\mathcal{V}_K$, for the global cluster mean field term $\check{x}^K$ by  \begin{equation*}
	\bar{x}^{K,q}\triangleq \mathrm{vec}(\bar{x}^{K,q}_1,\bar{x}^{K,q}_2, . . .,\bar{x}^{K,q}_K), 
\end{equation*}
where $\bar{x}^{K,q}_p$ denotes the estimation of  the local cluster mean field term $\check{x}^K_p$ by the agents in cluster $\mathcal{C}_q$.

 For agents in cluster $\mathcal{C}_q$, if $p\in\mathcal{N}_q$, then agents can obtain the local cluster mean field term $\check{x}^{K}_p$ by communication, i.e, 
\begin{equation}\label{e072701}
	 \bar{x}^{K,q}_p=\check{x}^{K}_p,
\end{equation}
else agents should make their own local estimation by mean field approximations according to \eqref{e010303}, i.e, 
\begin{equation}\label{e031701}
	\begin{split}
		&d\bar{x}^{K,q}_{p}(t)=[\tilde{A}_q(t)\bar{x}^{K,q}_p(t)+\tilde{G}_q(t)\bar{x}^{K,q}(t)]dt,\quad \bar{x}_p^K(0)=\bar{m}_p.
	\end{split}
\end{equation}
In order to represent the dynamics of the cluster mean field estimator more compact, the following notation is introduced for $q,p\in\mathcal{V}_K$:
\begin{align*}
	E_{qp}&\triangleq 1_n\otimes e_{ij},\\
	E^K_q&\triangleq \mathrm{rows}(E_{q1},E_{q2}, . . .,E_{qK}),\\
	\bar{E}_{qp}&\triangleq 1_n-E_{qp},\\ \bar{E}^K_q&\triangleq \mathrm{rows}(\bar{E}_{q1},\bar{E}_{q2}, . . .,\bar{E}_{qK}),\\
	Z^K&\triangleq \mathrm{diag}(B_1R_1^{-1}B_1^T\bar{K}_1,. . .,B_KR_K^{-1}B_K^T\bar{K}_K).
\end{align*}
By \eqref{e072701} and \eqref{e031701}, the dynamics of the designed cluster mean field estimator for agent $\mathcal{A}_i,i\in\mathcal{I}_q,q\in\mathcal{V}_K$, is given by 
\begin{align}
			d\bar{x}^{K,q}(t)=&E^K_q\odot d\check{x}^K(t) + \bar{E}^K_q\odot  [\tilde{A}^K(t) + \tilde{G}^K(t)]\bar{x}^{K,q}(t)dt, \quad  \bar{x}^{K,q}(0)=  E^K_q\odot x^K_{0}+\bar{E}^K_q\odot \bar{m}^K.\label{e031702}
\end{align} 


\subsection{Distributed controllers}
Based on the above discussion, we propose the following distributed controller for agent $\mathcal{A}_i,i\in\mathcal{I}_q,q\in\mathcal{V}_K$:
\begin{align}
	\check{u}_i(t)
	=&-R_q^{-1}B_q^TP_q(t)\check{x}_i(t)-R_q^{-1}B_q^T\bar{K}_q(t)\bar{x}^{K,q}(t),\label{e012109}
\end{align}
where $\bar{x}^{K,q}$ is the estimated global cluster mean field term of cluster $\mathcal{C}_q$ given by the cluster mean field estimator \eqref{e031702}. Under this distributed controller, the average distributed control of cluster $\mathcal{C}_q$ is given by
\begin{align}
	\check{u}_q^K(t)
	=& - R_q^{-1}B_q^TP_q(t)\check{x}_q^K(t) - R_q^{-1}B_q^T\bar{K}_q(t)\bar{x}^{K,q}(t).\label{e012108} 
\end{align}
Let  $$\bar{x}^K\triangleq\mathrm{vec}(\bar{x}^{K,1},\bar{x}^{K,2}, . . .,\bar{x}^{K,K})$$ be the estimated global cluster mean field term of all clusters in the large-population system. Therefore the closed-loop systems under the distributed controller are given by 
\begin{subequations}\label{e012301}
	\begin{align}	
		\begin{split}
			d\check{x}_i(t)=&[\tilde{A}_q(t)\check{x}_i(t)+\bar{G}_q(t)\check{x}^K(t)-B_qR_q^{-1}B_q^T\bar{K}_q(t)\bar{x}^{K,q}(t) ]dt+\Sigma_q dw_i(t),\quad \check{x}_i(0)=x_{i0},
		\end{split}\\
	\begin{split}\label{e012102}
			d\check{x}_q^K(t)=&[\tilde{A}_q(t)\check{x}_q^K(t)+\bar{G}_q\check{x}^K(t)- B_qR_q^{-1}B_q^T\bar{K}_q(t)\bar{x}^{K,q}(t) ]dt + \Sigma_q dw_q^K(t),\quad \check{x}_q^K(0)=x^K_{q0},
	\end{split}\\ 
			\begin{split}\label{e031901}
				d\check{x}^K(t)=&\big\{[\tilde{A}^K(t)+G^K]\check{x}^K(t)-Z^K(t)\bar{x}^{K}(t)\big\}dt+\Sigma^Kdw^K(t),\quad \check{x}^K(0)=x^K_0.
			\end{split}
	\end{align}
\end{subequations}


\subsection{Asymptotic optimality}
Next we will show the asymptotic optimality property of the designed distributed controller. We first define two kinds of estimation errors of the designed estimator \eqref{e031702} with respect to the realized cluster mean field terms under the centralized controller \eqref{e122909} and distributed controller \eqref{e012109}, respectively. For $q,p\in\mathcal{V}_K$, denote
\begin{align*}
\check{\xi}^{K,q}_p&\triangleq \check{x}^K_p-\bar{x}^{K,q}_{p},\\ \hat{\xi}^{K,q}_p&\triangleq \hat{x}^K_p-\bar{x}^{K,q}_p,\\
\check{\xi}^{K,q}&\triangleq \mathrm{vec}(\check{\xi}^{K,q}_{1},\check{\xi}^{K,q}_{2}, . . .,\check{\xi}^{K,q}_{K}),\\  
\hat{\xi}^{K,q}&\triangleq \mathrm{vec}(\hat{\xi}^{K,q}_{1},\hat{\xi}^{K,q}_{2}, . . .,\hat{\xi}^{K,q}_{K}),\\
\check{\xi}^K&\triangleq \mathrm{vec}(\check{\xi}^{K,1},\check{\xi}^{K,2}, . . .,\check{\xi}^{K,K}),\\
\hat{\xi}^K&\triangleq\mathrm{vec} (\hat{\xi}^{K,1},\hat{\xi}^{K,2}, . . .,\hat{\xi}^{K,K}).
\end{align*}

To prove the asymptotic optimality of the distributed controller \eqref{e012109}, we first obtain the following lemma, which shows that the estimator is asymptotically unbiased with respect to the number of agents.
\begin{lemma} \label{lemma1}
Assume (A1) holds. Let $C_1=\min_{q\in\mathcal{V}_K}N_q$. The cluster mean field estimator \eqref{e031701} is asymptotically unbiased, i.e.,
\begin{equation}\label{e032101}
	\mathbb{E}|
		\check{\xi}^{K}_q(t)|^2=O(\frac{1}{C_1}),\quad \mathbb{E}|
		\hat{\xi}^{K}_q(t)|^2=O(\frac{1}{C_1}),\quad q\in\mathcal{V}_K,
\end{equation}
and $\bar{x}^K$ and $\check{x}^K$ is bounded in the mean-square sense, i.e.,
\begin{equation}\label{e041701}
	\mathbb{E}|
		\check{x}^{K}(t)|^2=O(1),\quad 	\mathbb{E}|
		\bar{x}^{K}(t)|^2=O(1),\quad q\in\mathcal{V}_K.
\end{equation}

\end{lemma}

The following two lemmas are used to decompose the social cost function.

\begin{lemma}\label{lemma2}
	 Let $\zeta_i=x_i-x^K_q$, $v_i=u_i-u^K_q$, and $\mathbf{v}=\mathrm{vec}(v_1,v_2, . . .,v_N)$. Then the social cost \eqref{e122702} can be rewritten as 
	\begin{align}
		&J_{\rm soc}^{\rm (N)}(\mathbf{u})=J^1(u^K)+J^2(\mathbf{v}), \label{e012110}
	\end{align}
	where
	\begin{subequations}
	\begin{align}
	\begin{split}
			J^1(u^K)\triangleq&\sum_{q\in\mathcal{V}_K}\sum_{i\in\mathcal{I}_q}\mathbb{E}\Big\{\int_0^T \big[|x^K_q(t)-\bar{\Gamma}_qx^K(t)|^2_{Q_{q}} +|u^K_q(t)|^2_{R_{q}}\big]dt + |x^K_q(T) - \bar{\Gamma}_qx^K(T)|^2_{H_{q}}\Big\},
	\end{split}\\
	\begin{split}
			J^2(\mathbf{v})\triangleq &\sum_{q\in\mathcal{V}_K}\sum_{i\in\mathcal{I}_q}\mathbb{E}\int_0^T \big[|\zeta_i(t)|^2_{Q_q}+|v_i(t)|^2_{R_q}\big]dt\label{e032106}+\sum_{q\in\mathcal{V}_K}\sum_{i\in\mathcal{I}_q}\mathbb{E}|\zeta_i(T)|^2_{H_{q}}.
	\end{split}
	\end{align}
		\end{subequations}

\end{lemma}

\begin{lemma}\label{lemma3}
	Let	$\hat{\zeta}_i=\hat{x}_i-\hat{x}^K_q$, $\hat{v}_i= \hat{u}_i-\hat{u}^K_q$, $\check{\zeta}_i=\check{x}_i-\check{x}^K_q$ and $\check{v}_i= \check{u}_i-\check{u}^K_q$. The following equations hold for $i\in\mathcal{I}_q,q\in\mathcal{V}_K$:
	\begin{equation}\label{e032104}
		\hat{\zeta}_i(t)=\check{\zeta}_i(t),\quad \hat{v}_i(t)=\check{v}_i(t), 	
	\end{equation}
and 
 \begin{equation}\label{e032105}
 	J^2(\hat{\mathbf{v}})=J^2(\check{\mathbf{v}}).	
 \end{equation}

\end{lemma}

By Lemma \ref{lemma1}, \ref{lemma2} and \ref{lemma3}, the asymptotic optimality property of the designed distributed controller is given as follows. 
 \begin{theorem}\label{theorem3}
 	 Assume (A1) holds. The distributed controller \eqref{e012109} has the property of  asymptotically  social optimality, i,e.,
 	\begin{align}
 	\frac{J_{\rm soc}^{\rm (N)}(\mathbf{\check{u}})}{N}	-\frac{J_{\rm soc}^{\rm (N)}(\mathbf{\hat{u}})}{N}=O(\frac{1}{\sqrt{C_1}}).
 	\end{align}
 \end{theorem}

\section{Conclusion}\label{sec3}
This paper studies the problem of mean field social control in the case of heterogeneous agents following the direct approach. A set of asymptotically optimal distributed controllers is designed by constructing a cluster mean field estimator for each agent. In this paper we consider the finite cluster case with additive noise. An interesting work in the future is to further consider the heterogeneous mean field system with multiplicative noise in the infinite cluster case by graphon theory.

\end{document}